\documentclass[12pt,twoside]{article}
\usepackage{amsmath,amsbsy,amsfonts}
\newtheorem{theorem}{Theorem}[section]
\newtheorem{lemma}[theorem]{Lemma}

\newtheorem{corollary}[theorem]{Corollary}
\newtheorem{remark}[theorem]{Remark}

\newcommand{\qed}{\hfill\rule{2mm}{2mm}}

\title{Existence of least energy nodal solution for a  Schr\"odinger-Poisson system in bounded domains\footnote{Partially supported by INCT-MAT, Casadinho/PROCAD 552464/2011-2}}

\author{Claudianor O. Alves\footnote{C.O. Alves was partially supported by CNPq/Brazil 303080/2009-4, coalves@dme.ufcg.edu.br}
\,\,\, and \,\,\,  Marco A.S. Souto\footnote{M.A.S. Souto was
supported by CNPq/Brazil 304652/2011-3, marco@dme.ufcg.edu.br} \\
Universidade Federal de Campina Grande\\
Unidade Acad\^emica de Matem\'atica e Estat\'istica\\
CEP:58429-900, Campina Grande - PB, Brazil.}

\date{}

\begin{document}
\maketitle



\begin{abstract}

We prove the existence of least energy nodal solution for a class of Schr\"odinger-Poisson  system in a bounded domain $\Omega \subset \mathbb{R}^3$ with nonlinearity having  a subcritical growth.

\end{abstract}

\section{Introduction}

This paper was motivated by some works that have appeared in recent
years concerning  with the nonlinear Schr\"{o}dinger-Poisson
system
$$
\left\{ \begin{array}{l}
 -i\frac{\partial\psi}{\partial t} = - \Delta \psi +
\phi(x)\psi  - |\psi|^{p-2}\psi  \mbox{ in $\Omega$},\\
-\Delta \phi=|\psi|^2  \mbox{ in $\Omega$},\\ \phi= \psi=0 \mbox{ on
} \partial \Omega,
\end{array}
\right.\eqno{(NSP)}
$$
where $\Omega \subset \mathbb R^3$ is a bounded domain with smooth boundary, $2<p< 2^*=6$ and
$\psi : \overline{\Omega} \to \mathbb{C}$ and $\phi: \overline{\Omega} \to \mathbb{R}$ are unknown functions.

The first equation in $(NSP)$, called Schr\"odinger equation,
describes quantum (non-relativistic) particles interacting with the
eletromagnetic field generated by the motion. An interesting
Schr\"odinger equation class is when the potential $\phi(x)$ is
determined by the charge of wave function itself, that is, when the
second equation in $(NSP)$ (Poisson equation) holds.

Knowledge of the solutions for the elliptic equation
$$
\left\{ \begin{array}{ll}
  - \Delta u +
 \phi u = f(u)  \mbox{ in $\Omega$},\\
-\Delta \phi=u^2  \mbox{ in $\Omega$}, \\
u, \phi=0, \,\,\, \partial \Omega
\end{array}
\right.\eqno{(SP)}
$$
 has a great importance in the study of stationary solutions
 $\psi(x,t) = e^{-it}u(x)$ of $(NSP)$ and it contains two kinds of nonlinearities: the first one
is $\phi(x)u$ and concerns the interaction with the electric field.
This term is nonlocal, since the electrostatic potential $\phi(x)$
depends also on the wave function. The second nonlinearity is
$f(u)$. For more information involving physical situations where
$(SP)$ appears, we cite the papers of Benci-Fortunato \cite{bf},
Bokanowski \& Mauser \cite{BM}, Mauser \cite{Mauser}, Ruiz
\cite{Ruiz}, Ambrosetti-Ruiz \cite{AM} and S'anchez \& Soler
\cite{Sanchez}.

An important fact involving system $(SP)$ is that this class of system can be transformed into a Schr\"{o}dinger equation  with a nonlocal term (see, for instance, \cite{ap,G,Ruiz,zz}), which  allows to use variational methods. Effectively, by the Lax-Milgram Theorem, given $u\in H_0^1(\Omega)$, there exists
a unique $\phi=\phi_u \in H_0^1(\Omega)$ such that
$$
-\Delta \phi=u^2.
$$
By using standard arguments, we have that $\phi_u$ verifies the following properties (for a proof see
\cite{Daprile1,Ruiz, zz}):
\begin{lemma}{\label{lm1}} For any $u\in H^1_0(\Omega)$, we have
\begin{itemize}
\item[$i)$]  there exists $C>0$ such that
$||\phi_u||\leq C||u||^2$ and
$$
\int_{\Omega}|\nabla \phi_u|^2dx=\int_{\Omega}\phi_u u^2dx\leq
C ||u||^4 \quad \forall\, u\in H^1_0(\Omega);
$$
where $||u||^2=\int_{\Omega}|\nabla u|^2dx$.

\item[$ii)$] $\phi_u\geq 0$ $\forall u\in H^1_0(\Omega)$;

\item[$iii)$] $\phi_{tu}=t^2\phi_u$, $\forall t>0$ and $u\in H^1_0(\Omega)$;

\item[$iv)$] if $u_n \rightharpoonup u$ in $H_0^1(\Omega)$, then
$\phi_{u_n} \rightharpoonup \phi_u$ in $H^1_0(\Omega)$ and
$$
\lim_{n \rightarrow +\infty} \int_{\Omega}\phi_{ u_n}  u_n ^2dx=
\int_{\Omega}\phi_u u^2dx.
$$
\end{itemize}
\end{lemma}

Therefore, $(u,\phi) \in H^1_0(\Omega)\times H^1_0(\Omega)$ is a solution of $(SP)$ if, and
only if, $\phi = \phi_u$ and $u \in H^1_0(\Omega)$ is a weak solution of the nonlocal problem
$$
\left\{ \begin{array}{l}
  - \Delta u + \phi_u u = f(u)  \mbox{ in }\Omega,\\
u=0, \mbox { on }\partial \Omega .
\end{array}
\right.\eqno{(P)}
$$

Now, we would like to mention that the existence of solutions for problem $(P)$ can be made via variational methods, because if the nonlinearity $f$ belongs to $C^1(\mathbb R,\mathbb R)$ and satisfies
\begin{enumerate}
\item[$(f_1)$\ ]
$\displaystyle\lim_{s\rightarrow 0 } \frac{f(s)}{s} = 0$;
\item[$(f_2)$\ ] $\displaystyle\lim_{|s|\rightarrow +\infty } \frac{f(s)}{s^5} = 0$,
\end{enumerate}
the Lemma \ref{lm1} gives that the functional $J:
H^1_0(\Omega)\rightarrow \mathbb R$  given by
$$
J(u)=\frac 12 ||u||^2 +\frac 14\int_{\Omega}\phi_u u^2dx-
\int_{\Omega}F(u)dx,
$$
where
$$
F(s)=\int_0^s f(t)dt,
$$
belongs to $C^1(H^{1}_{0}(\Omega),\mathbb{R})$ and
\[
J'(u)v =  \int_{\Omega}\nabla u \nabla v dx + \int_{\Omega}\phi_u
uvdx - \int_{\Omega}f(u)vdx \,\,\,\, \forall u,v \in H^1_0(\Omega).
\]
Hence, critical points of $J$ are the weak solutions for nonlocal problem $(P)$.

From  the above commentaries, we have that system $(SP)$ has a
nontrivial solution if, and only if, $(P)$ has a nontrivial
solution. This way, in the last years, many authors that studied the
system $(SP)$ have focused their attention on problem $(P)$ aiming
to establish existence and nonexistence of solutions, multiplicity
of solutions, ground state solutions, radial and nonradial
solutions, semiclassical limit and concentrations of solution for
the case where $\Omega = \mathbb{R}^{N}$, see the papers of
Azzollini \& Pomponio \cite{ap},  Cerami \& Vaira \cite{Cerami},
Coclite \cite{C}, D'Aprile \& Mugnai \cite{Daprile1,Daprile2},
d'Avenia \cite{Davenia}, Ianni \cite{Ianni}, Kikuchi \cite{K1}, and
Zhao \& Zhao \cite{zz}. For the case where $\Omega$ is a bounded
domain, we would like to cite the papers of Siciliano \cite{G}, Ruiz
\& Siciliano \cite{Ruiz-Sic} and Pisani \& Siciliano \cite{LG}. In
all those papers, the solutions found are nonnegative. However,
related to nodal ( or sign-changing ) solution, we found few papers,
see for example,  Ianni \cite{Ianni2} and Kim \& Seok \cite{Kim}. In
\cite{Ianni2} and \cite{Kim} the existence of nodal solutions have
been established at balls centered  origin  or in whole
$\mathbb{R}^{3}$.

Motivated by papers above, we are interested in finding nodal solution for system $(SP)$, by assuming only that $\Omega \subset \mathbb{R}^{3}$ is a bounded domain with smooth boundary. Once that we will apply variational methods and term  $\int_{\Omega}\phi_u u^2dx$ is homogeneous of degree
4, the corresponding Ambro\-se\-tti-Rabinowitz condition on $f$ is
the following:
\begin{itemize}
\item[(AR)] There exists $\theta >4$ such that
\[
0 < \theta F(s) \leq sf(s)\quad \forall  s\in \mathbb{R}\setminus \{0\}.
\]
\end{itemize}
This condition is important not only to ensure that the functional
$J$ has the mountain pass geometry, but also to guarantee that the
Palais-Smale, or  Cerami, sequences associated with $J$ are bounded.
We recall that (AR) implies a weaker condition: there exist
$\theta>4$ and $C_1,C_2>0$ such that
\begin{equation}\label{ar}
F(s) \geq C_1|s|^\theta - C_2,\quad \forall\, s\in \mathbb{R}.
\end{equation}
However, we consider here another much weaker one, namely,
\begin{itemize}
\item[$(f_3)$\ ] $\displaystyle\lim_{s\rightarrow +\infty} \frac{F(s)}{s^4}
=+\infty$.
\end{itemize}
Moreover, we also assume that the nonlinearity $f$ satisfies
\begin{itemize}
\item[$(f_4)$\ ] $\displaystyle \frac{f(s)}{s^3} $ is increasing in $|s|>0$.
\end{itemize}

\vspace{0.5 cm}

\begin{remark}\label{rem} The condition $(f_4)$ implies that $H(s)=sf(s)-4F(s)$ is a non-negative function, increasing increasing in $|s|$ with
$$
sH'(s)=s^{2}f'(s)-3f(s)s>0 \,\,\, \mbox{for any} \,\, |s|>0.
$$
\end{remark}

\newpage
Our main result is the following

\begin{theorem}{\label{fth}}
Suppose that   $f$  satisfies $(f_1)- (f_4)$. Then problem ($P$)
possesses at least energy nodal solution, which has precisely two
nodal domains.
\end{theorem}

In the proof of Theorem \ref{fth}, we prove that functional $J$ assumes
a minimum value on the nodal set
$$
\mathcal M=\{u\in \mathcal N: J'(u)u^+=J'(u)u^-=0 \mbox { and }
u^{\pm}\neq 0 \}
$$
where $u^{+}=\max\{u(x),0\}$, $u^{-}(x)=\min\{u(x),0\}$ and
$$
\mathcal{N} = \{ u\in
H^1_0(\Omega)\setminus\{0\}\, :\, J'(u)u=0\}.
$$
More precisely, we prove that there is $w \in \mathcal{M}$ such that
$$
J(w)=\inf_{u \in \mathcal{M}}J(u).
$$
After, motivated by argument used in Bartsch, Weth \& Willem
\cite{BWW}, we use a deformation lemma to prove that $w$ is a
critical point of $J$, and so, $w$ is a least energy nodal solution
for $(SP)$ with exactly two nodal domains.

Since  $J$ has the
nonlocal term $\int_{\Omega}\phi_u u^2dx$, if $u$ is a nodal
solution for $J$, we have that
$$
J'(u^+)u^+ =-\int_{\Omega}\phi_{u^{-}}(u^{+})^{2} <0\,\,\, \mbox{and} \,\,\, J'(u^-)u^- =-\int_{\Omega}\phi_{u^{+}}(u^{-})^{2} <0.
$$
From this, some arguments used to prove the existence of nodal solutions for problem like
$$
\left\{ \begin{array}{l}
  - \Delta u = f(u)  \mbox{ in }\Omega,\\
u=0, \mbox { on }\partial \Omega
\end{array}
\right.\eqno{(P_1)}
$$
can not be used, and so, a careful analysis  is necessary in a lot of estimates, see Section 2 for details.

Before to conclude this introduction, we would like to cite the papers of Alves \cite{A}, Alves \& Soares \cite{AS1,AS2}, Bartsch, Weth and Willem \cite{BWW}, Bartsch \& Weth \cite{BW}, Bartsch, Liu \& Weth \cite{BLW},  Castro, Cossio \& Neuberger \cite{ACN}, Zou \cite{Zou} and their references, where existence of nodal solution has been studied for problem related to $(P_1)$.

The paper is organized as follows. In Section 2, we show some estimates involving functions that change sign, with the most of them being new for problem $(P)$. The Section 3 is devoted to prove the main result Theorem \ref{fth}.

\section { Important estimates }

In what follows, we denote by $\mathcal{N}$  the Nehari manifold associated with $J$, that is,
$$
\mathcal{N} = \{ u\in
H^1_0(\Omega)\setminus\{0\}\, :\, J'(u)u=0\}.
$$
A critical point $u_0\neq 0$ of $J$ is a ground state of $(P)$ if
$$
J(u_0) = \inf_\mathcal{N} J(u).
$$
Since we are looking for least energy nodal solutions (or
sign-changing solutions), our goal is to prove the existence of a
critical point for $J$ in the set
$$
\mathcal M=\{u\in \mathcal N: J'(u)u^+=J'(u)u^-=0 \mbox { and }
u^{\pm}\neq 0 \}.
$$

Let us start with some technical lemmas.

\begin{lemma}\label{lema2}
There exists $\rho>0$ such that
\begin{enumerate}
\item [(i)] $J(u)\geq ||u||^2/4$ and $||u||\geq \rho, \forall u\in \mathcal N$;
\item [(ii)] $||w^\pm||\geq\rho, \,\, \forall w \in \mathcal M$.
\end{enumerate}
\end{lemma}

\noindent {\bf{Proof:}} From $(f_4)$ and Remark \ref{rem}, for any $u\in \mathcal N$
$$
4J(u)=4J(u)-J'(u)u= ||u||^2+ \int_{\Omega}[uf(u)-4F(u)]dx\geq
||u||^2
$$
and so,
$$
J(u)\geq ||u||^2/4 \,\,\, \forall u \in \mathcal N .
$$

From $(f_1)$ and $(f_2)$, there is $C>0$ such that
$$
f(s)s\leq \frac {\lambda_1}2 s^2+Cs^6, \mbox { for all } s\in
\mathbb R.
$$
where $\lambda_1$ is the first eigenvalue of $(-\Delta , H_{0}^{1}(\Omega)).$ Since $J'(u)u=0$,
$$
||u||^2<||u||^2+\int_\Omega \phi_{u}|u|^2dx=\int_{\Omega}uf(u)dx\leq
\frac {\lambda_1}2\int_{\Omega}u^2dx +C \int_{\Omega}u^6dx.
$$
Then by Sobolev embeddings,
$$
||u||^2< \frac{1}{2} ||u||^2+\hat C ||u||^6,
$$
from where it follows that
$$
 ||u||\geq \rho  \,\,\, \forall u \in \mathcal N,
$$
where $\rho = \left(\frac 1{2\hat C}\right)^{\frac{1}{4}}$, finishing the proof of $(i).$

If $w \in \mathcal M$, we have that $J'(w)w^{\pm}=0$. Then, a simple computation gives $J'(w^{\pm})w^{\pm}<0$, which implies
$$
||w^{\pm}||^2<||w^{\pm}||^2+\int_\Omega \phi_{w^{\pm}}(w^{\pm})^2dx<\int_\Omega
f(w^{\pm})w^{\pm}dx.
$$
As in the item $(i)$, we can deduce that $||w^{\pm}||\geq \rho$.

 \qed

 \begin{lemma}\label{lema2x}
If $(w_n)$ is a bounded sequence in  $\mathcal M$ and $p\in
(2,6)$, we have
$$
\liminf_n \int_{\Omega}|w^\pm_n|^{p}dx>0.
$$
\end{lemma}

\noindent {\bf{Proof:}} From $(f_1)$ and $(f_2)$, given
$\varepsilon>0$ there exists $C>0$ such that
$$
f(s)s\leq \varepsilon\lambda_1 s^2+C|s|^p+\varepsilon s^6, \mbox {
for all } s\in \mathbb R.
$$
Since $w_n\in \mathcal M $, by Lemma \ref{lema2}
$$
\rho^{2} \leq||w^{\pm}_n||^2<\int_{\Omega}w^+_nf(w^+_n)dx \leq
\varepsilon\lambda_1\int_{\Omega}(w^+_n)^2dx
+C\int_{\Omega}|w^+_n|^pdx+\varepsilon \int_{\Omega}(w^+_n)^6dx
$$
that is,
$$
\rho^2 \leq \varepsilon \left( \lambda_1\int_{\Omega}(w^{\pm}_n)^2dx + \int_{\Omega}(w^{\pm}_n)^6dx  \right) + C\int_{\Omega}|w^{\pm}_n|^pdx.
$$
Using the boundedness of $(w_n)$, there is $C_1$ such that
$$
\rho^2 \leq \varepsilon C_1 + C\int_{\Omega}|w^{\pm}_n|^pdx.
$$
Fixing $\varepsilon = \frac{\rho^2}{2C_1}$, we get
$$
\int_{\Omega}|w^{\pm}_n|^pdx \geq \frac{\rho^2}{2C},
$$
showing that
$$
\liminf_{n} \int_{\Omega}|w^{\pm}_n|^pdx \geq \frac{\rho^2}{2C} >0.
$$
\qed

\begin{lemma}\label{lema3}
Let  $v\in H_0^1(\Omega)$ with $v^{\pm}\neq 0$. Then,
there are $t,s>0$ such that $J'(tv^++sv^-)v^+=0$ and
$J'(tv^++sv^-)v^-=0$.
\end{lemma}

\noindent {\bf{Proof:}} It what follows, we consider the vector field
$$
V(s,t)=\left (J'(tv^++sv^-)(tv^+),J'(tv^++sv^-)(sv^-)\right ).
$$
from $(f_1)-(f_3)$, a  straightforward computation yields that there are $0<r<R$ such that
$$
J'(rv^++sv^-)(rv^+), \,\,\, J'(tv^++rv^-)(rv^-)>0, \,\,\, \forall s,t\in[r,R]
$$
and
$$
J'(Rv^++sv^-)(Rv^+), \,\,\, J'(tv^++Rv^-)(Rv^-)<0, \,\,\, \forall s,t\in[r,R].
$$
Now, the lemma follows applying  Miranda theorem \cite{Miranda}. \qed

\vspace{0.5 cm}

Hereafter, for $v\in H_0^1(\Omega)$ with  $v^{\pm}\neq 0$, we consider  the functions $h^v:[0,+\infty)\times [0,+\infty) \to \mathbb R$ given by
$$
h^{v}(t,s)=J(tv^++sv^-)
$$
and $\Phi^v:[0,+\infty)\times [0,+\infty) \to \mathbb R^2$ defined as
$$
\Phi^v(t,s)=\left (\frac{\partial h^v}{\partial t}(t,s),\frac{\partial
h^v}{\partial s}(t,s)\right )=\left
(J'(tv^++sv^-)v^+,J'(tv^++sv^-)v^-\right ).
$$
Since $f$ is a $C^1$ function, it follows that $\Phi^v $ is also a $C^1$ map. Moreover, it is easy to check that if $(t,s)$ is a critical point of $h^v$, then
$$
\begin{array}{ll}
h^v(t,s)= &  h^v(t,s)-\frac 14 \left\langle \nabla h^v(t,s),(t,s) \right\rangle \\
\mbox{} & \mbox{}\\
\mbox{} & =\frac 14t^2||v^+||^2+ \displaystyle \frac 14\int_\Omega [f(tv^+)tv^+-4F(tv^+)]dx+ \\
\mbox{} & \mbox{}\\
\mbox{} & \,\,\,\,\,\,\, \frac 14s^2||v^-||^2+ \displaystyle \frac 14\int_\Omega [f(sv^-)tv^--4F(sv^-)]dx.
\end{array}
$$
\vspace{0.5 cm}
\begin{lemma}\label{lema4}
If $w\in \mathcal M$,
\begin{enumerate}
\item [(a)] $h^w(t,s)<h^w(1,1)=J(w)$, for all $s,t\geq 0$ such that
$(s,t)\neq (1,1)$;
\item [(b)] $\det (\Phi^w)'(1,1)>0$.
\end{enumerate}
\end{lemma}

\noindent {\bf{Proof:}} Once that $w\in \mathcal M$,  we have
$J'(w)w^+=J'(w)w^-=0$, and so,
$$
||w^+||^2+\int_\Omega \phi_{w^+}(w^+)^2dx+\int_\Omega
\phi_{w^-}(w^+)^2dx=\int_\Omega f(w^+)w^+dx
$$
and
$$
||w^-||^2+\int_\Omega \phi_{w^-}(w^-)^2dx+\int_\Omega
\phi_{w^+}(w^-)^2dx=\int_\Omega f(w^-)w^-dx.
$$
These equalities imply that $(1,1)$ is a critical point of $h^w$. On the other hand, condition $(f_3)$ leads to  the limit
$$
\lim_{|(t,s)| \to \infty }h^w(t,s)=-\infty,
$$
which implies $h^w$ assumes a global maximum in some $(a,b)$.

First of all, we claim that $a,b>0$. If $b=0$,  we have
$$
J(aw^+)\geq J(tw^+), \,\,\, \forall t>0
$$
and then $J'(aw^+)(aw^+)=0$, or equivalently,
$$
a^2||w^+||^2+a^4\int_\Omega \phi_{w^+}(w^+)^2dx=\int_\Omega
f(aw^+)aw^+dx.
$$
Since $J'(w^+)w^+=J'(w)w^+-\int_\Omega \phi_{w^-}(w^+)^2dx<0$, we derive
$$
||w^+||^2+\int_\Omega \phi_{w^+}(w^+)^2dx<\int_\Omega f(w^+)w^+dx
$$
and so,
$$
\left (1-\frac 1{a^2} \right )||w^+||^2<\int_\Omega \left ( \frac
{f(w^+)w^+}{(w^+)^4}-\frac {f(aw^+)aw^+}{(aw^+)^4} \right
)(w^+)^4dx.
$$
If $a>1$ the left side in this inequality is positive while, from
$(f_4)$, the right side is negative. This information gives that $ a\leq 1$. Now, combining the Remark \ref{rem} with the fact that $a \leq 1$, we get
\begin{eqnarray*}
h^w(a,0)=J(aw^+)= J(aw^+)-\frac 14J'(aw^+)(aw^+)=\\=\frac 14
a^2||w^+||^2+\frac 14\int_\Omega [f(aw^+)aw^+-4F(aw^+)]dx\\\leq\frac
14 ||w^+||^2+\frac 14\int_\Omega [f(w^+)w^+-4F(w^+)]dx\\<\frac 14
||w^+||^2+\frac 14\int_\Omega [f(w^+)w^+-4F(w^+)]dx+\\+\frac 14
||w^-||^2+\frac 14\int_\Omega [f(w^-)w^- -4F(w^-)]dx=\\=J(w)-\frac
14J'(w)w=J(w)=h(1,1)
\end{eqnarray*}
that is,
$$
h^w(a,0) < h^w(1,1)
$$
which is absurd, because $(a,0)$ is a global maximum point for $h^w$. The same type of argument works to show that $a\not=0$, and the proof of claim is done.

The second claim is $0<a,b\leq 1$. In fact, since $(a,b)$ is another
critical point of $h^w$,
$$
a^2||w^+||^2+a^4\int_\Omega \phi_{w^+}(w^+)^2dx+a^2b^2\int_\Omega
\phi_{w^-}(w^+)^2dx=\int_\Omega f(aw^+)aw^+dx
$$
and
$$
b^2||w^-||^2+b^4\int_\Omega \phi_{w^-}(w^-)^2dx+a^2b^2\int_\Omega
\phi_{w^+}(w^-)^2dx=\int_\Omega f(bw^-)bw^-dx.
$$
Without loss of generality, we will suppose that $a\geq b$. From this,
$$
a^2||w^+||^2+a^4\int_\Omega \phi_{w^+}(w^+)^2dx+a^4\int_\Omega
\phi_{w^-}(w^+)^2dx\geq\int_\Omega f(aw^+)aw^+dx
$$
leading to
$$
\left( \frac 1{a^2}-1\right )||w^+||^2\geq\int_\Omega \left ( \frac
{f(aw^+)aw^+}{(aw^+)^4}-\frac {f(w^+)w^+}{(w^+)^4} \right
)(w^+)^4dx.
$$
If $a>1$ the left side in this inequality is negative, but from
$(f_4)$, the right side is positive, thus we can deduce that $a \leq 1$.

To conclude the proof of item $(a)$,  we will show that $h^w$ does not
have global maximum  in $[0,1]\times [0,1]\setminus \{(1,1)\}$. From definition of $h^w$,
\begin{eqnarray*}
h^w(a,b)=\frac 14a^2||w^+||^2+\frac 14\int_\Omega
[f(aw^+)aw^+-4F(aw^+)]dx+\\ \frac 14b^2||w^-||^2+\frac 14\int_\Omega
[f(bw^-)bw^--4F(bw^-)]dx.
\end{eqnarray*}
Then, if $0<a,b\leq 1$ and $(a,b)\neq (1,1)$,
\begin{eqnarray*}
h^w(a,b)< \frac 14||w^+||^2+\frac 14\int_\Omega
[f(w^+)w^+-4F(w^+)]dx+\\ \frac 14||w^-||^2+\frac 14\int_\Omega
[f(w^-)w^--4F(w^-)]dx= h^w(1,1)
\end{eqnarray*}
showing that,
$$
h^w(a,b) < h^w(1,1)
$$
and thereby, the proof of item $(a)$ is complete.

The proof of item $(b)$ is the following. By a simple calculation
$$
\det (\Phi^w)'(1,1)=G(w^+)G(w^-)-4\left [\int_\Omega
\phi_{w^-}(w^+)^2dx\right ]^2
$$
where
$$
G(v)=\int_\Omega [f'(v)v^2-f(v)v]dx-2\int_\Omega \phi_{v}v^2dx.
$$
From Remark \ref{rem}
$$
G(v)\geq 2\left [\int_\Omega f(v)vdx-\int_\Omega \phi_{v}v^2dx\right
].
$$
Once that
$$
\int_\Omega f(w^+)w^+dx-\int_\Omega
\phi_{w^+}(w^+)^2dx=||w^+||^2+\int_\Omega \phi_{w^-}(w^+)^2dx
$$
and
$$
\int_\Omega \phi_{w^-}(w^+)^2dx=\int_\Omega \phi_{w^+}(w^-)^2dx,
$$
we have that
$$
G(w^+)> 2 \int_{\Omega} \phi_{w^-}(w^+)^2dx
$$
and
$$
G(w^-)> 2 \int_{\Omega} \phi_{w^-}(w^+)^2dx.
$$
Combining the above informations, it follows that $\det (\Phi^w)'(1,1)>0$. \qed

\begin{corollary} \label{C1}  Let $v \in H^{1}_0(\Omega)$ be a function verifying
$$
v^\pm \not =0 \,\,\, \mbox{and} \,\,\, J'(v)v^{\pm} \leq 0 .
$$
Then, there are $t,s \in [0,1]$ such that
$$
tv^+ + s v^- \in \mathcal{M}.
$$
\end{corollary}
\noindent {\bf Proof.} An immediate consequence of the arguments used in the proof of Lemma \ref{lema4}. \qed

\section { Existence of least energy nodal solution.  }

In this section, our main goal is to prove the Theorem \ref{fth}. In what follows, we denote by $c_0$ the infimum of $J$ on $\mathcal M$, that
is,
$$
c_0=\inf_{v\in \mathcal M} J(v).
$$
From Lemma \ref{lema2}(i), we deduce that $c_0>0$.

Let $(w_n)$ be a sequence in $\mathcal M$ such that
$$
\lim_{n} J(w_n)=c_0.
$$
Still from Lemma \ref{lema2}(i), $(w_n)$ is a bounded
sequence. Hence, without loss of generality, we can suppose that there is $w \in H_0^1(\Omega)$ verifying
$$
w_n \rightharpoonup w \,\,\, \mbox{in} \,\,\, H_0^{1}(\Omega),
$$

$$
w_n \to w \,\,\, \mbox{in} \,\,\,  L^p(\Omega) \,\,\ \forall \, p \in [1,2^{*})
$$
and
$$
w_n(x) \to w(x) \,\,\, \mbox{a.e. in} \,\, \Omega.
$$

The condition $(f_2)$ combined with the \emph{compactness lemma of Strauss}
\cite[Theorem A.I, p.338]{bl} gives
$$
\lim_n \int_{\Omega}|w_n^{\pm}|^{p}dx = \int_{\Omega}|w^{\pm}|^{p}dx,
$$
$$
\lim_n \int_\Omega w_n^\pm f(w_n^\pm)dx=\int_\Omega w^\pm f(w^\pm)dx
$$
and
$$
\lim_n \int_\Omega  F(w_n^\pm)dx=\int_\Omega F(w^\pm)dx,
$$
from where it follows together with Lemma \ref{lema2x} that $w^\pm\neq 0$. Then, by Lemma \ref{lema3} there are $t,s>0$ verifying
$$
J'(tw^++sw^-)w^+=0 \,\,\, \mbox{and} \,\,\, J'(tw^++sw^-)w^-=0.
$$
Next, we will show that $t,s\leq 1$. Since
$J'(w_n)w_n^\pm=0$,
$$
||w^+_n||^2+\int_\Omega \phi_{w^+_n}(w^+_n)^2dx+\int_\Omega
\phi_{w^-_n}(w^+_n)^2dx=\int_\Omega f(w^+_n)w^+_ndx
$$
and
$$
||w^-_n||^2+\int_\Omega \phi_{w^-_n}(w^-_n)^2dx+\int_\Omega
\phi_{w^+_n}(w^-_n)^2dx=\int_\Omega f(w^-_n)w^-_ndx.
$$
Taking the limit in the above equalities, we obtain
$$
||w^+||^2+\int_\Omega \phi_{w^+}(w^+)^2dx+\int_\Omega
\phi_{w^-}(w^+)^2dx\leq\int_\Omega f(w^+)w^+dx
$$
and
$$
||w^-||^2+\int_\Omega \phi_{w^-}(w^-)^2dx+\int_\Omega
\phi_{w^-}(w^+)^2dx\leq\int_\Omega f(w^-)w^-dx.
$$
Once that
$$
J'(tw^++sw^-)(tw^+)=J'(tw^++sw^-)(sw^-)=0,
$$
it follows that
$$
t^2||w^+||^2+t^4\int_\Omega \phi_{w^+}(w^+)^2dx+t^2s^2\int_\Omega
\phi_{w^-}(w^+)^2dx=\int_\Omega f(tw^+)tw^+dx
$$
and
$$
s^2||w^-||^2+s^4\int_\Omega \phi_{w^-}(w^-)^2dx+t^2s^2\int_\Omega
\phi_{w^-}(w^+)^2dx=\int_\Omega f(sw^-)sw^-dx.
$$
Now, without loss of generality, we will suppose that $s\geq t$.
Under this condition,
$$
s^2||w^-||^2+s^4\int_\Omega \phi_{w^-}(w^-)^2dx+s^4\int_\Omega
\phi_{w^-}(w^+)^2dx\geq\int_\Omega f(sw^-)sw^-dx
$$
and then
$$
\left( \frac 1{s^2}-1\right )||w^-||^2\geq\int_\Omega \left ( \frac
{f(sw^-)sw^-}{(sw^-)^4}-\frac {f(w^-)w^-}{(w^-)^4} \right
)(w^-)^4dx.
$$
If $s>1$, the left side in this inequality is negative, but from
$(f_4)$, the right side is positive, thus we must have $s \leq 1$, which also implies that $t \leq 1$.

Our next step is show that $J(tw^++sw^-)=c_0$.  Recalling that \linebreak $tw^++sw^-\in
\mathcal M$, we derive that
$$
c_0 \leq J(tw^++sw^-)=J(tw^++sw^-)-\frac{1}{4}
J'(tw^++sw^-)(tw^++sw^-).
$$
Thus,
$$
c_0 \leq \left (J(tw^+)-\frac 14
J'(tw^+)(tw^+)\right )+\left (J(sw^-)-\frac 14
J'(sw^-)(sw^-)\right ).
$$
From Remark \ref{rem},
$$
J(tw^+)-\frac 14 J'(tw^+)(tw^+)  \leq J(w^+)-\frac 14 J'(w^+)(w^+)
$$
and
$$
J(sw^-)-\frac 14 J'(sw^-)(sw^-) \leq J(w^-)-\frac 14 J'(w^-)(w^-).
$$
Hence,
$$
c_0 \leq \left (J(w^+)-\frac 14
J'(w^+)(w^+)\right )+\left (J(w^-)-\frac 14
J'(w^-)(w^-)\right ).
$$
Using Fatous' Lemma combined again with Remark \ref{rem},
$$
c_0 \leq  J(tw^+ + sw^-) \leq \liminf_{n} \left (J(w_n) -\frac{1}{4}J'(w_n)w_n \right)=\lim_n J(w_n)=c_0
$$
from where it follows that
$$
c_0 = J(tw^+ + sw^-).
$$

Until this moment, we have proved that there exists a $w_o=tw^+ + sw^-\in
\mathcal M$, such that $J(w_o)=c_0$. In what follows, let us denote $w_o$ by $w$, consequently
$$
J(w)=c_0 \,\,\, \mbox{and} \,\,\, w\in \mathcal M.
$$

To conclude the proof of Theorem \ref{fth}, we claim that $w$ is a critical point for
functional $J$. If it is not true, there exist $\alpha >0$ and  $v_0 \in H^1_0(\Omega)$ with  $||v_0||=1$
satisfying
$$
J'(w)v_0=2\alpha>0.
$$
Since $J'$ is continuous, we fix $r>0$ such that
$$
J'(v)v_0>\alpha, v^\pm \neq 0, \mbox{ for all } \,\, v \in B_r(w) \subset H_0^{1}(\Omega).
$$
From now on, fix $D=(\xi,\chi)\times (\xi,\chi) \subset \mathbb
R^2$ with  $0<\xi<1<\chi$ such that
\begin{enumerate}
\item [(i)] $(1,1)\in D$ and $\Phi^w(t,s)=0$ in $\overline D$ if, and only if,
$t=s=1$;
\item [(ii)] $c_0 \notin h^w(\partial D)$;
\item [(iii)] $\{tw^++sw^-: (t,s)\in \overline D\} \subset B_r(w)$;
\end{enumerate}
where $h^w$ and $\Phi^w$ were defined in Lemma \ref{lema4}. Since $J$ is
continuous, we can fix $r'>0$ such that
$$
\mathcal B=\overline
{B_{r'}(w)}\subset B_r(w)
$$
and
$$
\mathcal B\cap \{tw^++sw^-: (t,s)\in
\partial D\}=\emptyset.
$$
Consider the continuous mapping $\rho:H_0^1(\Omega) \to
[0,+\infty)$, defined by
$$
\rho(u)=\mbox{dist}(u, \mathcal B^c).
$$
Moreover, set the bounded Lipschitz vector field
$V:H_0^1(\Omega)\to H_0^1(\Omega)$ given by
$$
V(u)=-\rho(u)v_0.
$$
For each $u\in H_0^1(\Omega)$, we denote by $\eta(\tau)=\eta(\tau,u)$ the
unique solution of ODE
\[
\left \{\begin{array}{l}
\eta'(\tau)=V(\eta(\tau)), \, t>0\\
\eta(0)=u.
\end{array}\right .
\]
Observe that
\begin{enumerate}
\item [(1) ] if $u\notin \mathcal B$, $\eta(\tau,u)=u$, for all $t$;
\item [(2)] if $u\in \mathcal B$,  $\tau \mapsto J(\eta(\tau,u))$ is decreasing and $\eta(\tau,u)\in \mathcal B$, for all $\tau>0$;
\item [(3)] there exists $\tau_o>0$ such that $J(\eta(\tau,w))\leq J(w)-((r'\alpha)/2))\tau$, for all $0\leq \tau \leq \tau_o$.
\end{enumerate}

The item $(1)$ is an immediate consequence from  the definition of $\rho$. The item $(2)$ follows from the inequality
$$
J'(\eta(\tau))\eta'(\tau) \leq -\rho(\eta(\tau))\alpha<0, \,\,\, \forall \eta(\tau)\in \mathcal B.
$$
To verify (3), fix $\tau_o>0$ such that
\[
||\eta(\tau,w)-w||\leq \frac {r'} 2, \mbox { for all } |\tau|\leq
\tau_o.
\]
Thus,
\[
\frac d{dt}J(\eta(\tau,w)) \leq  -\rho(\eta(\tau))\alpha\leq -\frac
{r'\alpha}2.
\]
Integrating in $[0,\tau_0]$, we have
\[
J(\eta(\tau_0,w))\leq J(w)-\frac {r'\alpha}2\tau_0.
\]
Now, consider $\gamma:\overline D\to H_0^1(\Omega)$ given by.
$$
\gamma(t,s)=\eta(\tau_o,tw^++sw^-).
$$
It is easy to see that
$$
\max_{(t,s)\in\overline D} J( \gamma(t,s)) <c_0,
$$
because
$$
J(\gamma(t,s))\leq h^w(t,s)<c_0 \,\,\, \forall (t,s)\in \overline
D \setminus \{(1,1)\}
$$
and
$$
J(\gamma(1,1))\leq
J(w)-((r'\alpha)/2))\tau_o<c_o.
$$
Consequently $\gamma(\overline D)\cap\mathcal M=\emptyset$.

On the other hand, setting $\Psi:\overline D \to \mathbb R^2$ by
$$
\Psi(t,s)=(t^{-1}J'(\gamma(t,s))(\gamma(t,s)^+),s^{-1}J'(\gamma(t,s))(\gamma(t,s)^-)),
$$
we derive that
$$
\Psi(t,s)=\left
(J'(tw^++sw^-)w^+,J'(tw^++sw^-)w^-\right )=\Phi^w(t,s) \,\,\, \forall (t,s)\in \partial D.
$$
Then, using the Brouwer's topological degree
$$
d(\Psi,D,(0,0))=d(\Phi^w,D,(0,0))=\mbox{ sgn}(\det (\Phi^w)'(1,1))=1
$$
which yields $\Psi$ has a zero $(a,b)$ in $D$. Thereby, there is $ (a,b) \in D$ verifying
$$
J'(\gamma(a,b))(\gamma(a,b)^{\pm})=0,
$$
that is, $\gamma(a,b)\in \mathcal M$ which is a contradiction. From this, $w$ is a critical point of $J$, and so, a nodal solution for problem $(P)$. Now, we will show that $w$ has exactly two nodal domains, to this end, we assume by contradiction that
$$
w=u_1+u_2+u_3 \,\,\,
$$
with
$$
u_i \not= 0, u_1 \geq 0, u_2 \leq 0  \,\,\, \mbox{and} \,\,\ suppt(u_i) \cap suppt(u_j) = \emptyset \,\, i \not= j \,\, (i,j=1,2,3).
$$
Setting $v = u_1 + u_2$, we see that $v^{\pm} \not= 0$. Moreover, using the fact that $J'(w)=0$, it follows that
$$
J'(v)(v^{\pm}) \leq 0.
$$
By Corollary \ref{C1}, there are $t,s \in (0,1]$ such that
$$
tv^+ + sv^{-} \in \mathcal{M}
$$
or equivalently,
$$
tu_1 + su_2 \in \mathcal{M},
$$
and so,
$$
J(tu_1+su_2) \geq c_0.
$$
On the other hand, repeating the same type of argument explored in the proof of Lemma \ref{lema4} combined with the fact that $u_3 \not=0$, we find
$$
J(tu_1 + su_2) < J(w)=c_0,
$$
obtaining a contradiction. This way, $u_3=0$, and $w$ has exactly two nodal domains. \qed

\vspace{0.5 cm}

\noindent \textbf{Acknowledgments.} { The authors are  grateful to the referees
for a number of helpful comments for improvement in this article.}


\begin{thebibliography}{99}

\bibitem{A} C.O. Alves, Multiplicity of multi-bump type nodal solutions for a class of elliptic problems in $\mathbb{R}^{N}$. Top. Meth. Nonlinear Anal. 34 (2009), 231-250.

\bibitem{AS1} C.O. Alves, S.H.M. Soares, On the location and profile of spike-layer
nodal solutions to nonlinear Schr\"odinger equations. J. Math. Anal. Appl. 296 (2004), 563 - 577.

\bibitem{AS2} C.O. Alves, S.H.M. Soares, Nodal solutions for singularly perturbed
equations with critical exponential growth, J. Differential Equations 234 (2007), 464 - 484.

\bibitem{AM} {A. Ambrosetti} \& {R. Ruiz, }  {Multiple bound states for the {S}chr\"odinger-{P}oisson
              problem}, {Commun. Contemp. Math.} {10} {(2008)}
  {391--404}.


   \bibitem{ap}
    {A. Azzollini} \& {A. Pomponio, }
   {Ground state solutions for the nonlinear
              {S}chr\"odinger-{M}axwell equations},
{J. Math. Anal. Appl.}
     {345}
{(2008)}
  {90--108}.


\bibitem{BWW}{T. Bartsch, T. Weth} \& { M. Willem,} {Partial symmetry of least energy nodal solution to some variational problems}, {Journal D'Analyse Math\'ematique}{ 1 (2005) 1-18}


\bibitem{BW} T. Bartsch \& T. Weth, Three nodal solutions of singularly perturbed elliptic equations
on domains without topology, Ann. Inst. H. Poincar\'e Anal. Non
Lin\'eaire 22 (2005), 259�-281.

\bibitem{BLW} T. Bartsch, Z. Liu \&  T. Weth, Sign changing solutions of superlinear  {S}chr\"odinger
equations, Comm. Partial Diff. Equations 29 (2004), 25�-42.










    \bibitem{bf}
    {V. Benci } \& {D. Fortunato, }
 {An eigenvalue problem for the Schr\"odinger-Maxwell equations},
{Top. Meth. Nonlinear Anal.} {11} {(1998)}
     {283--293}.


\bibitem{bl} {H. Berestycki} \& {P.L. Lions, }
    {Nonlinear scalar field equations, I - existence of a ground state },
    {Arch. Rat. Mech. Analysis},
   {82},
     {(1983)},
   {313--346}.


\bibitem{BM} {O. Bokanowski} \& {N.J. Mauser,}{ Local approximation of the Hartree-Fock exchange potential: a deformation approach, ${\text M}^3$AS 9 (1999), 941-961.}





\bibitem{ACN} A. Castro, J. Cossio \& J. Neuberger, A sign-changing solution for a
superlinear Dirichlet problem, Rocky Mountain Journal of Mathematics.
27, 4 (1997), 1041- 1053.



\bibitem{Cerami}
 {G. Cerami} \& {G. Vaira, }
 {Positive solutions for some non-autonomous
              {S}chr\"odinger-{P}oisson systems},
{J. Differential Equations}
     {248} {(2010)}
  {521--543}.


\bibitem{C}
 {G.M. Coclite, }
 {A multiplicity result for the nonlinear
              {S}chr\"odinger-{M}axwell equations},
{Commun. Appl. Anal.}
     {7}
{(2003)}
  {417--423}.





\bibitem{Daprile1}
 {T. D'Aprile} \& {D. Mugnai},
 {Solitary waves for nonlinear {K}lein-{G}ordon-{M}axwell and
              {S}chr\"odinger-{M}axwell equations},
{Proc. Roy. Soc. Edinburgh Sect. A} {134} {(2004)}
  {893--906}.





\bibitem{Daprile2}
 {T. D'Aprile} \& {D. Mugnai},
 {Non-existence results for the coupled  {K}lein-{G}ordon-{M}axwell  equations},
{Adv. Nonlinear Stud.} {4} {(2004)}
  {307--322}.





\bibitem{Davenia}
 {P. d'Avenia, }
 {Non-radially symmetric solutions of nonlinear {S}chr\"odinger
              equation coupled with {M}axwell equations},
{Adv. Nonlinear Stud.}  {(2002)}
 {2}
  {177--192}.




     \bibitem{G}
    {G. Siciliano, }
 {Multiple positive solutions for a  Schr\"odinger-Poisson-Slater system},
{J. Math. Analysis and Appl.} {365} {(2010)}
     {288--299}.

\bibitem{K1}
 {H. Kikuchi},
 {On the existence of a solution for elliptic system related to
              the {M}axwell-{S}chr\"odinger equations},
{Nonlinear Anal.}     {67} {(2007)} {1445--1456}.



\bibitem{Ianni}
 {I. Ianni} \& {G. Vaira},
 {On concentration of positive bound states for the
              {S}chr\"odinger-{P}oisson problem with potentials},
{Adv. Nonlinear Stud.}     {8} {(2008)}  {573-595}.

\bibitem{Ianni2}
 {I. Ianni}, {Sign-Changing radial solutions for the Schr\"odinger-Poisson-Slater problem},
{arXiv:1108.2803v1 } .



\bibitem{Kim}
 {S. Kim } \& {J. Seok},
 {On nodal solutions of the Nonlinear Schr\"odinger-Poisson equations},
{Comm. Cont. Math.}     {14} {(2012)}  {12450041-12450057}.


\bibitem{Miranda} {C. Miranda}, {Un' osservazione su un teorema di Brouwer,}{ Bol. Un. Mat. Ital., 3 (1940)} 5-7.



\bibitem{Mauser} {N.J. Mauser,}  {The Schr\"odinger-Poisson-$X_\alpha$ equation,} {Applied Math. Letters 14 (2001)},
759-763.

\bibitem{LG} {L. Pisani} \& {G. Siciliano,} { Note on a Schr\"odinger-Poisson system in a bounded domain,
Appl. Math. Lett.} {21} {(2008)}  {521-528} .




     \bibitem{Ruiz}
    {D. Ruiz, }
 {The Schr\"odinger-Poisson equation under the effect of a nonlinear local term },
{J. Funct. Analysis} {237} {(2006)}  {655--674}.



\bibitem{Ruiz-Sic}
 {D. Ruiz} \& {G. Siciliano,}
 {A note on the Schr\"odinger-Poisson-{S}later equation on
              bounded domains}, {Adv. Nonlinear Stud.}  {8}
{(2008)}  {179--190}.





\bibitem{Sanchez}
{O. S\'anchez} \& {J. Soler}, {Long-time dynamics of the Schr\"odinger-Poisson-Slater system,
J. Statistical Physics 114 (2004), 179-204.}


\bibitem{zz}
    {F. Zhao }\& {L. Zhao,}
 {Positive solutions for Schr\"{o}dinger-Poisson equations with a critical exponent},
{Nonlinear Anal.} {70} {(2009)}    {2150--2164}.



\bibitem{Zou} {W. Zou},\,{Sign-Changing critical point theory}, Springer, 2008.


\end{thebibliography}
\end{document}